\documentclass[12pt]{article}
\usepackage{amssymb}
\usepackage{epsf,graphicx,color}
\usepackage{amsmath}
\usepackage{tikz}
\usepackage[pdftoolbar=true,
            pdfmenubar=true,
            pdfpagemode=UseOutlines,
            bookmarksnumbered=true,
            linktocpage=true,
            colorlinks=false,
]{hyperref}
\begin{document}


\renewcommand{\theequation}{\mbox{\arabic{equation}}}

\newtheorem{assumption}{Assumption}
\newtheorem{theorem}{Theorem}
\newtheorem{lemma}{Lemma}
\newtheorem{definition}{Definition}
\newtheorem{corollary}{Corollary}
\newtheorem{proposition}{Proposition}
\newtheorem{remark}{Remark}
\newtheorem{conjecture}{Conjecture}
\newtheorem{example}{Example}

\renewcommand{\theassumption}{\mbox{A\arabic{assumption}}}
\renewcommand{\thetheorem}{\mbox{\arabic{section}.\arabic{theorem}}}
\renewcommand{\thelemma}{\mbox{\arabic{section}.\arabic{lemma}}}
\renewcommand{\thedefinition}{\mbox{\arabic{section}.\arabic{definition}}}
\renewcommand{\thecorollary}{\mbox{\arabic{section}.\arabic{corollary}}}
\renewcommand{\theproposition}{\mbox{\arabic{section}.\arabic{proposition}}}
\renewcommand{\theremark}{\mbox{\arabic{section}.\arabic{remark}}}
\renewcommand{\theconjecture}{\mbox{\arabic{section}.\arabic{conjecture}}}
\renewcommand{\theexample}{\mbox{\arabic{section}.\arabic{example}}}

\newcommand{\proof}{{\bf Proof}. }

\newcommand{\diag}{{\rm diag}}
\newcommand{\rank}{{\rm rank}}

\newcommand{\qed}{~\hfill$\square$}



\def\R{\mathbb{R}}
\def\H{\mathbb{H}}
\def\S{\mathbb{S}}
\def\Rquer{\overline{\mathbb{R}}}
\def\N{\mathbb{N}}
\def\Z{\mathbb{Z}}
\def\C{\mathbb{C}}
\def\Q{\mathbb{Q}}
\def\P{\mathcal{P}}


\def\AA{\mathfrak{A}}
\def\BB{\mathfrak{B}}
\def\BBquer{\overline{\mathfrak{B}}}
\def\DD{\mathfrak{D}}
\def\CC{\mathit{C}}
\def\EE{\mathfrak{E}}
\def\E{\mathcal{E}}
\def\FF{\mathfrak{F}}
\def\GG{\mathfrak{G}}
\def\II{\mathfrak{I}}
\def\LL{\mathfrak{L}}
\def\PP{\mathfrak{P}}
\def\pp{\mathfrak{p}}
\def\RR{\mathfrak{R}}
\def\NN{\mathcal{N}}
\def\KK{\mathcal{K}}
\def\HH{\mathcal{H}}
\def\TT{\mathcal{T}}


\def\Po{{\cal P}o}


\def\dann{\mbox{ }\Longrightarrow\mbox{ }}
\def\Aequi{\Longleftrightarrow}
\def\PPfeil{\Longrightarrow}
\def\Pfeil{\longrightarrow}
\def\pfeil{\rightarrow}


\def\eps{\varepsilon}
\def\sig{\sigma}
\def\theta{\vartheta}
\def\rho{\varrho}
\def\om{\omega}
\def\Om{\Omega}


\def\kom{\mbox{, }}
\def\dif{\,\mbox{d}}
\def\pkt{{\scriptscriptstyle \bullet}}
\def\ov{\overline}
\def\leer{\emptyset}
\def\ohne{\backslash}
\def\Eins{{1\hspace{-0.39em}1}}                             
\def\esup{ ess\,sup\,}
\def\einf{ ess\,inf\,}

\def\Frac#1#2{{\displaystyle \frac{#1}{#2}}}
\def\sfrac#1#2{{\textstyle \frac{#1}{#2}}}
\def\Binom#1#2{{\displaystyle \binom{#1}{#2}}}

\def\norm#1#2#3{{ {\parallel {#1} \parallel}_{#2}^{#3}}}
\def\Norm#1{ \parallel {#1} \parallel}


\def\Int#1#2{{\displaystyle \int_{#1}^{#2} }}
\def\Sum#1#2{{\displaystyle \sum_{#1}^{#2} }}
\def\Prod#1#2{{\displaystyle \prod_{#1}^{#2} }}
\def\Bigcap#1#2{{\displaystyle \bigcap_{#1}^{#2} }}
\def\Bigcup#1#2{{\displaystyle \bigcup_{#1}^{#2} }}
\def\Bigtimes#1#2{{\displaystyle \overset{#2}{\underset{#1}{\text{\raisebox{-2pt}{\huge \mbox{}\hspace{2pt}$\times$} }}} }}
\def\Bigotimes#1#2{{\displaystyle \bigotimes_{#1}^{#2} }}

\def\Lim#1{{\displaystyle \lim_{#1}}}
\def\Limn{{\displaystyle \lim_{n \pfeil \infty}}}

\def\p-Limn{{\displaystyle \mathcal{P-}\lim_{n \pfeil \infty}\, }}

\def\Limm{{\displaystyle \lim_{m \pfeil \infty}}}
\def\Limk{{\displaystyle \lim_{k \pfeil \infty}}}
\def\Limt{{\displaystyle \lim_{t \pfeil \infty}}}
\def\Limh{{\displaystyle \lim_{h \pfeil 0}}}
\def\Liminfn{{\displaystyle \liminf_{n \pfeil \infty}}}
\def\Limsupn{{\displaystyle \limsup_{n \pfeil \infty}}}
\def\Limsupk{{\displaystyle \limsup_{k \pfeil \infty}}}

\def\Max#1{{\displaystyle \max_{#1}}}
\def\Min#1{{\displaystyle \min_{#1}}}
\def\Sup#1{{\displaystyle \sup_{#1}}}
\def\Inf#1{{\displaystyle \inf_{#1}}}

\def\ds{\displaystyle}

\renewcommand{\dagger}{*}

\title{\textbf{Topological approach to mathematical programs with switching constraints}}

\author{
V. Shikhman
\thanks{
Department of Mathematics, Chemnitz University of Technology,
Reichenhainer Str. 41, 09126
Chemnitz, Germany; e-mail:  vladimir.shikhman@mathematik.tu-chemnitz.de (corresponding author).
 }
}

\date{}

\maketitle

\begin{abstract}
\noindent
We study mathematical programs with switching constraints (MPSC) from the topological perspective. Two basic theorems from Morse theory are proved. Outside the W-stationary point set, continuous deformation of lower level sets can be performed. However, when passing a W-stationary level, the topology of the lower level set changes via the attachment of a $w$-dimensional cell. The dimension $w$ equals the W-index of the nondegenerate W-stationary point. The W-index depends on both the number of negative eigenvalues of the restricted Lagrangian's Hessian and the number of bi-active switching constraints. As a consequence, we show the mountain pass theorem for MPSC.
Additionally, we address the question if the assumption on the nondegeneracy of W-stationary points is too restrictive in the context of MPSC. It turns out that all W-stationary points are generically nondegenerate. Besides, we examine the gap between nondegeneracy and strong stability of W-stationary points. A complete characterization of strong stability for W-stationary points by means of first and second order information of the MPSC defining functions under linear independence constraint qualification is provided. In particular, all  bi-active  Lagrange  multipliers of a strongly stable W-stationary point cannot  vanish.


\end{abstract}

\setcounter{lemma}{0}
\setcounter{theorem}{0}
\setcounter{definition}{0}
\setcounter{example}{0}
\setcounter{proposition}{0}
\setcounter{corollary}{0}

\section{Introduction}

We consider the following mathematical program with switching constraints:
\begin {equation}
  \label{eq:math-program}
   \mbox{MPSC}[f,h,g,F_1,F_2]:\quad \min \,\, f(x) \quad \mbox{ s.t.} \quad x\in M[h,g,F_1,F_2]
\end{equation}
with
    \[
       M[h,g,F_1,F_2]=\left\{x \in \R^n \left|
       \begin{array}{l}
                           h_i(x) =0, i \in I, g_j(x) \geq 0, j \in J, \\
                          F_{1,m}(x) \cdot F_{2,m}(x) = 0, m=1,\ldots, k 
                           \end{array}\right. \right\},
  \]
  where $f \in C^2(\R^n,\R)$, $ h \in C^2(\R^n, \R^{|I|})$, $ g \in C^2(\R^n, \R^{|J|})$,
  $F_1, F_2 \in C^2(\R^n,\R^k)$. For simplicity, we write $M$ for $M[h,g,F_1,F_2]$ if no confusion is possible. MPSC has been introduced in \cite{mehlitz:2019} by arguing that switching structures, which demand at most one control to be active at any time instance, appear rather frequently in optimal control. Discretization of such optimal control problems naturally leads to MPSC. In \cite{mehlitz:2019}, various stationarity concepts for MPSC were suggested, such as W-, M-, and S-stationarity. Under suitable MPSC-tailored constraint qualifications, necessary optimality conditions were derived. Numerical aspects of MPSC were handled in \cite{kanzow:2019}, where different relaxation schemes for switching constraints are discussed and compared.

%
%

The goal of this paper is the study of MPSC from the topological point of view as it has been pioneered by H. Th. Jongen for nonlinear programming in \cite{jongen:1977} and popularized in \cite{jongen:2000}. Subsequent studies in this direction are summarized in \cite{shikhman:2012}, where, in particular, mathematical programs with complementarity constraints (MPCC) and mathematical programs with vanishing constraints (MPVC) are considered. 
%
%
  The main question here is how the topological properties of lower level sets
\[
   M^a=\{ x \in M \,|\, f(x) \leq a\}
\]
change as the level $a \in \R$  varies. For MPSC, it turns out that the concept of W-stationarity is adequate in order to describe these changes.
In particular, within this context, we present two basic theorems from {\bf Morse theory}, cf. \cite{jongen:2000, milnor:1963}.
First, we show that, for $a < b$, the set $M^a$ is a strong deformation retract of $M^b$ if the set 
\[
   M^b_a=\{ x \in M \,|\, a \leq f(x) \leq b\}
\]
does not contain W-stationary points, see Theorem \ref{theorem:main}(a). 
As a consequence, the homotopy type of the lower level sets $M^a$ and $M^b$ are equal. 
Second, if $M^b_a$ contains exactly one nondegenerate W-stationary point,
then $M^b$ is shown to be homotopy-equivalent to $M^a$ with a $w$-cell attached, see Theorem \ref{theorem:main}(b).
Here, the dimension $w$ is the so-called W-index. It depends on both the number of negative eigenvalues of the restricted Lagrangian's Hessian and the number of bi-active switching constraints. The latter fact constitutes the main difference to the cases where feasible set is described by complementarity or vanishing constraints. We remind that for MPCC the C-index of topologically relevant C-stationary points depends on the number of negative pairs of bi-active Lagrange multipliers, see \cite{jongen:2009}, \cite{ralph:2011}. The same is true for MPVC, where just negative pairs of bi-active Lagrange multipliers of T-stationary points matter, see \cite{dorsch:2012}.  For switching constraints however {\bf all pairs of bi-active Lagrange multipliers contribute to W-index}.
A global interpretation of the deformation and cell-attachment theorems is as follows. Suppose that the feasible set is compact
and  connected, that the linear independence constraint qualification (LICQ) holds, and that all W-stationary points are nondegenerate with pairwise different functional values. Then, passing a level corresponding to a local minimizer, a connected component of the lower level set is created. Different components can only be connected by attaching one-dimensional cells. This shows the existence of at least $(r-1)$ W-stationary saddle points with W-index equal to one, where $r$ is the number of local minimizers, see also \cite{floudas_jongen:2005, jongen:2000}. This mountain pass result for MPSC is shown in Theorem \ref{thm:mrel}.

%
%
We point out that the crucial assumption in Morse theory is the nondegeneracy of W-stationary points. A W-stationary point is called nondegenerate, see Definition \ref{def:nondeg_w_stationary}, if (ND1) LICQ is satisfied,
  (ND2) the Lagrange multipliers corresponding to active inequality constraints are positive, (ND3) the Lagrange multipliers corresponding to bi-active switching constraints do not vanish, and (ND4) the restricted Hessian of the Lagrangian is nonsingular.
Are the requirements ND1-ND4 too restrictive to be assumed for W-stationary points of MPSC? First, we show that all W-stationary points of a generic MPSC are nondegenerate, see Theorem \ref{thm:ND-gen}. Genericity means that the set of MPSC defining functions $(f,g,h,F_1,F_2)$, for which all W-stationary points are nondegenerate, is open and -dense with respect to the so-called strong (or Whitney) $C^2$-topology, see e.\,g. \cite{hirsch:1976}. Second, we examine the gap between nondegeneracy and strong stability of W-stationary points. Strongly stable W-stationary points  in the sense of M. Kojima \cite{kojima:1980} do not only remain locally unique under sufficiently small $C^2$-perturbations of the MPSC defining functions $(f, h, g, F_1, F_2)$, but can also be  continuously tracked back. All nondegenerate W-stationary points are strongly stable, but not vice versa. In this paper, we present a full {\bf characterization of strong stability} for W-stationary points by means of first and second order information of the MPSC defining functions  under LICQ, see Theorem \ref{thm:char-ss}. 
In particular, ND3 is necessary for strong stability of W-stationary points, i.\,e. {\bf all bi-active Lagrange multipliers  cannot  vanish}. This new issue is in strong contrast e.\,g. with the characterization of strong stability for C-stationary points in MPCC, see \cite{jongen:2012}. In MPCC, strong stability includes cases where one of the bi-active Lagrange multipliers may vanish. We conclude that in absence of active inequality constraints, nondegeneracy of W-stationary points is equivalent to their strong stability. 

The article is organized as follows. In Section 2 we provide notation and auxiliary results which will be used later.
Section 3 contains the exposition of Morse theory for MPSC including the proofs of the deformation and cell-attaching theorem.
In Section 4 we fully characterize strongly stable W-stationary points of MPSC and draw comparisons to their nondegeneracy.

Our notation is standard. The $n$-dimensional Euclidean space is denoted by $\R^n$ with norm $\|\cdot\|$. We denote the set of nonnegative numbers by $\H$. The solution set of the basic switching relation is denoted by
\[
  \S = \{ \left(a, b\right) \in \R^2 \, | \, a \cdot b =0\}.
\]
$B(\bar x,r)$ stands for the ball $\left\{\left. x \in \R^n \,\right|\, \|x - \bar x\| \leq r\right\}$ for $\bar x \in \R^n$, $r > 0$.
Given a differentiable function $F:\R^n \longrightarrow \R^m$, $DF$ denotes its Jacobian matrix.
Given a differentiable function $f:\R^n \longrightarrow \R$, $Df$ denotes the row vector of partial derivatives 
and $D^T f$ stands for the transposed vector. 

\setcounter{lemma}{0}
\setcounter{theorem}{0}
\setcounter{definition}{0}
\setcounter{example}{0}
\setcounter{proposition}{0}
\setcounter{corollary}{0}

\section{Notations and auxiliary results}

Given $\bar x \in M$, we define the following index sets:
\[
  J_0(\bar x)=\{j \in J \,|\, g_j(\bar x) = 0\},
\]
\[
  \alpha(\bar x)=\{m \in \{1,\ldots k\} \,|\, F_{1,m}(\bar x) = 0, F_{2,m}(\bar x) \not = 0 \},
\]
\[
  \beta(\bar x)=\{m \in \{1,\ldots k\} \,|\, F_{1,m}(\bar x) = 0, F_{2,m}(\bar x) = 0 \},
\]
\[
  \gamma(\bar x)=\{m \in \{1,\ldots k\} \,|\, F_{1,m}(\bar x) \not = 0, F_{2,m}(\bar x) = 0 \}.
\]
The index set $J_0(\bar x)$ corresponds to the active inequality constraints and $\beta({\bar x})$ to the bi-active switching constraints at $\bar x$.

Without loss of generality, we assume throughout the whole article that
at the particular point of interest $\bar x \in M$ it holds:
\[
J_0(\bar x)=\{1, \ldots, |J_0(\bar x)|\}, \alpha(\bar x)=\{1, \ldots, |\alpha(\bar x)|\},
\]
\[
\gamma(\bar x)=\{|\alpha(\bar x)|+1, \ldots, |\alpha(\bar x)|+|\gamma(\bar x)|\}.
\]

We put $s= |I|+|\alpha(\bar x)|+|\gamma(\bar x)|$,
  $q= s + |J_0(\bar x)|$,
  $p=n - q - 2|\beta(\bar x)|$.

Let us start by stating the MPSC-tailored linear independence constraint qualification, which turns out to be the crucial assumption for all results to follow.  

\begin{definition}[LICQ, cf. \cite{mehlitz:2019}]
The linear independence constraint qualification (LICQ) for MPSC is said to hold at $\bar x \in M$ if the vectors
\[   \begin{array}{l}
        D h_i(\bar x), i \in I,
       D F_{1,m_\alpha}(\bar x), m_\alpha \in \alpha (\bar x),
       D F_{2,m_\gamma}(\bar x), m_\gamma \in \gamma (\bar x), \\
       D g_j(\bar x), j \in J_0(\bar x),
        D F_{1,m_\beta}(\bar x), D F_{2,m_\beta}(\bar x), m_\beta \in \beta (\bar x)
       
     \end{array}
  \]
are linearly independent.
\end{definition}

The assumption of LICQ is justified in the sense, that it generically holds on the MPSC feasible set. In order to formulate this assertion in mathematically precise terms, the space $C^2(\R^n,\R)$ will be topologized by means of the strong (or Whitney-) $C^2$-topology, denoted by $C^2_s$, cf. \cite{hirsch:1976, jongen:2000}. The $C^2_s$-topology is generated by allowing
perturbations of the functions and their derivatives up to second order which are controlled by means of continuous
positive functions. The product space $C^2(\R^n, \R^l) \cong C^2(\R^n, \R) \times \cdots \times C^2(\R^n, \R)$ will be topologized with the corresponding product topology.

\begin{theorem}[LICQ is generic]
\label{thm:licq-gen}
Let $\mathcal{F}$ denote the subset of $C^2(\R^n,\R^{|I|}) \times C^2(\R^n,\R^{|J|}) \times C^2(\R^n,\R^k) \times C^2(\R^n,\R^k)$ consisting of those MPSC defining functions for which LICQ holds at all feasible points. Then, $\mathcal{F}$ is $C^2_s$-open and -dense.
\end{theorem}

\proof Let us consider the following mathematical program with vanishing constraints:
\[
   \mbox{MPVC}[f,h,g,F_1,F_2] \quad \min \,\, f(x) \quad \mbox{s.t.} \quad  x\in N[h,g,F_1,F_2]
\]
with
  \[
       N[h,g,F_1,F_2]=\left\{x \in \R^n \left|
       \begin{array}{l}
                           h_i(x) =0, i \in I, g_j(x) \geq 0, j \in J, \\
                           F_{1,m}(x) \geq 0, F_{1,m}(x) \cdot F_{2,m}(x) \leq 0, m=1,\ldots, k 
                           \end{array}\right. \right\}.
  \]
Obviously, it holds $M[h,g,F_1,F_2] \subset N[h,g,F_1,F_2]$. The MPVC-tailored LICQ is shown to hold generically on the feasible set $N[h,g,F_1,F_2]$, see \cite[Theorem 2.1(i)]{dorsch:2012}. Moreover, the definitions of MPVC- and MPCC-tailored LICQ coincide on $M[h,g,F_1,F_2]$. The assertion thus follows. \qed

Among different stationarity concepts for MPSC discussed in \cite{mehlitz:2019} that of W-stationarity appears to be relevant for our topologically motivated studies. 

\begin{definition}[W-stationary point, cf. \cite{mehlitz:2019}]
   \label{def:w_stationary}
   A point $\bar x \in M$ is called W-stationary for MPSC if there exist real numbers (Lagrange multipliers)
   \[
   \bar \lambda_i, i \in I,
   \bar \mu_j, j \in J,
   \bar \sigma_{1, m}, \bar \sigma_{2, m}, m = 1, \ldots, k,
   \]
   such that:
    \begin{equation}
          \label{eq:w_crit_diff}
         \begin{array}{rcl}
               D f(\bar x) &=& \displaystyle \Sum{i \in I}{} \bar \lambda_i D h_i(\bar x)+\Sum{j \in J}{} \bar \mu_j D g_j(\bar x) \\
        &&\displaystyle +\Sum{m=1}{k} \left( \bar \sigma_{1,m} D F_{1,m}(\bar x)+
         \bar \sigma_{2,m} D F_{2,m}(\bar x) \right),
         \end{array}
    \end{equation}
    \begin{equation}
         \label{eq:w_crit_mu_active}
             \bar \mu_j \cdot g_j(\bar x)=0, j \in J,             
     \end{equation}
     \begin{equation}
         \label{eq:w_crit_mu_nonnegative}
                  \bar \mu_j \geq 0, j \in J,
     \end{equation}
     \begin{equation}
         \label{eq:w_crit_sigma_active}
             \bar \sigma_{j, m} \cdot F_{j, m}(\bar x) = 0, j=1,2, m=1,\ldots, k.
      \end{equation}
A vector $(\bar x, \bar \lambda, \bar \mu, \bar \sigma) \in M \times \R^{|I|} \times \R^{|J|} \times \R^{2k}$ satisfying (\ref{eq:w_crit_diff})-(\ref{eq:w_crit_sigma_active}) is called a W-stationary pair for MPSC.
\end{definition}

 In the case where LICQ holds at a W-stationary point $\bar x \in M$, the Lagrange multipliers in (\ref{eq:w_crit_diff}) are uniquely determined. It is well-known that under LICQ a local minimum of MPSC is W-stationary \cite{mehlitz:2019}. 

Let us now locally describe the MPSC feasible set under LICQ. This is done by an appropriate change of coordinates.

\begin{definition}[Coordinate system]
   \label{definition:local_coordinates}
  The feasible set $M$ admits a local $C^r$-coordinate system of $\R^n$ ($r \geq 1$) at $\bar x$
  by means of a $C^r$-diffeomorphism $\Phi:U \longrightarrow V$ with open $\R^n$-neighborhoods $U$ and $V$ of $\bar x$ and $0$,
  respectively, if it holds:
   \begin{itemize}
      \item[(i)]  $\Phi(\bar x)=0$,
      \item[(ii)] $\Phi(M \cap U)=
                  \left( \{0_{s}\} \times \H^{|J_0(\bar x)|}
                  \times \S^{|\beta(\bar x)|}
                  \times \R^{p} \right) \cap V$.
   \end{itemize}
\end{definition}

\begin{lemma}[Local structure]
   \label{lemma:feasible_set}
   Suppose that LICQ holds at $\bar x \in M$. Then $M$ admits a local $C^2$-coordinate system of $\R^n$ at $\bar x$.
\end{lemma}

\proof
Choose vectors $\xi_l \in \R^n$, $l= 1, \ldots, p$, which form - together with the vectors
\[   \begin{array}{l}
       D^T h_i(\bar x), i \in I,
       D^T F_{1,m_\alpha}(\bar x), m_\alpha \in \alpha (\bar x),
       D^T F_{2,m_\gamma}(\bar x), m_\gamma \in \gamma (\bar x), \\
       D^T g_j(\bar x), j \in J_0(\bar x),
        D^T F_{1,m_\beta}(\bar x), D^T F_{2,m_\beta}(\bar x), m_\beta \in \beta (\bar x)
     \end{array}
  \]
- a basis for $\R^n$.
Next we put

\begin{equation}
      \label{eq:standard_diffeo}
\left.        \begin{array}{lll}
y_i &=& h_i(x), i\in I, \\
y_{|I|+m_\alpha} &=& F_{1,m_\alpha}(x), m_\alpha \in \alpha(\bar x), \\
y_{|I|+m_\gamma} &=& F_{2,m_\gamma}(x), m_\gamma \in \gamma(\bar x),  \\
y_{s+j} &=& g_j(x), j \in J_0(\bar x), \\
y_{s+|J_0(\bar x)|+2 m_\beta -1} &=&
      F_{1,m_\beta}(x), \\
y_{s+|J_0(\bar x)|+2  m_\beta } &=&
      F_{2,m_\beta}(x), m_\beta =1, \ldots, |\beta(\bar x)|, \\
y_{n-p+l} &=& \xi^T_l( x- \bar x), l=  1, \ldots, p. \\
  \end{array}
\right\}
\end{equation}
or, shortly,
\begin{equation}
      \label{eq:standard_diffeo_1}
       y=\Phi(x).
\end{equation}
Note that $\Phi \in C^2(\R^n,\R^n)$, $\Phi(\bar x)=0$ and the Jacobian matrix $D\Phi(\bar x)$ is nonsingular in virtue of LICQ and the choice of $\xi_l$, $l= 1, \ldots, p$. By means of the implicit function theorem there exist open neighborhoods $U$ of $\bar x$ and $V$ of $0$ such that $\Phi:U \longrightarrow V$ is a $C^2$-diffeomorphism. By shrinking $U$, if necessary, we can guarantee that $J_0(x) \subset J_0(\bar x)$ and $\beta(x) \subset \beta(\bar x)$ for all $x \in M \cap U$. Thus, the property (ii) in Definition \ref{definition:local_coordinates} follows directly from the definition of $\Phi$. \qed

\begin{definition}
  \label{definition:standard_diffeo}
  We will refer to the $C^2$-diffeomorphism $\Phi$ defined by (\ref{eq:standard_diffeo}), (\ref{eq:standard_diffeo_1})
  as standard diffeomorphism.
\end{definition}

Let us interpret the first- and second-order derivatives of the objective function $f \circ \Phi^{-1}$ in new coordinates given by the standard diffeomorphism $\Phi$. For that, we define the Lagrange function
\[
        \begin{array}{rcl}
        L(x, \lambda, \mu, \sigma)&=&\displaystyle f(x)-
         \Sum{i \in I}{} \lambda_i h_i( x)-\Sum{j \in J}{} \mu_j g_j( x) \\
        &&\displaystyle -\Sum{m=1}{k} \left( \sigma_{1,m} F_{1,m}( x)+ \sigma_{2,m} F_{2,m}( x) \right).
        \end{array}
\]
At $\bar x \in M$ we further consider the local part of the feasible set:
\[
  \begin{array}{ll}
      M_0(\bar x)= \{ x \in \R^n \,\,| & h_i(x) =0, i \in I,
    g_j(x) =0, j \in J_0(\bar x),  \\ &
       F_{1,m_\alpha}(x) = 0, m_\alpha \in \alpha(\bar x), F_{2,m_\gamma}(x) = 0, m_\gamma \in \gamma(\bar x) , \\ &
    F_{1,m_\beta}(x) = 0, F_{2,m_\beta}(x) = 0, m_\beta \in \beta(\bar x)
       \}.
     \end{array}
\]
Obviously, $M_0(\bar x) \subset M$ and, in the case where LICQ holds at $\bar x$, $M_0(\bar x)$ is locally a
$p$-dimensional $C^2$-manifold. The tangent space of $M_0(\bar x)$ at $\bar x$ is
\[
    \begin{array}{ll}
      T_{\bar x}M_0(\bar x) = \{ \xi \in \R^n \,\,| &
        D h_i(\bar x) \, \xi=0, i \in I, D g_j(\bar x) \, \xi =0, j \in J_0(\bar x) \\ &
        D F_{1,m_\alpha}(\bar x) \, \xi = 0, m_\alpha \in \alpha(\bar x), \\ &
        D F_{2,m_\gamma}(\bar x) \, \xi= 0, m_\gamma \in \gamma(\bar x), 
     \\ &
    D F_{1,m_\beta}(\bar x) \, \xi = 0, D F_{2,m_\beta}(\bar x) \, \xi = 0, m_\beta \in \beta(\bar x)
         \}.
     \end{array}
\]

\begin{remark}[New coordinates]
 \label{remark:deriv_in_new_coordinates}
From the proof of Lemma \ref{lemma:feasible_set} it follows that the Lagrange multipliers $(\bar \lambda, \bar \mu, \bar \sigma)$ of a W-stationary point $\bar x$ are the corresponding partial derivatives of the objective function $f\circ\Phi^{-1}$ in new coordinates given by the standard diffeomorphism $\Phi$, cf. \cite[Lemma 2.2.1]{jongen:2004}. Moreover, its Hessian with respect to the last $p$ coordinates corresponds to the restriction of the Lagrange function's Hessian $D^2_{xx}L(\bar x, \bar \lambda,\bar \mu, \bar \sigma)$ on the respective tangent space $T_{\bar x}M_0(\bar x)$, cf. \cite[Lemma 2.2.10]{jongen:2004}.
\end{remark}

\section{Morse theory}

Our goal is to study topological changes of MPSC lower level sets as the level of the objective function varies. For the topological concepts we refer to \cite{jongen:2000, spanier:1966}. 
In order to provide intuition on deformation and cell-attachment in presence of switching constraints, we first consider the following simple Example \ref{ex:ca}.

\begin{figure}
    \centering
\begin{tikzpicture}
\draw[->] (-3,0) -- (xyz cs:x=-3,y=3);
\node at (-2.6,2.8){$x_2$};
\draw[->] (-4.5,1.5) -- (xyz cs:x=-1.5,y=1.5);
\node at (-1.7,1.2){$x_1$};
\draw[ultra thick] (-4.5,1.5) -- (xyz cs:x=-4,y=1.5);
\draw[ultra thick] (-3,0) -- (xyz cs:x=-3,y=0.5);
\node at (-3.75,0.75){$M^a$};

\draw[->] (1,0) -- (xyz cs:x=1,y=3);
\node at (1.4,2.8){$x_2$};
\draw[->] (-0.5,1.5) -- (xyz cs:x=2.5,y=1.5);
\node at (2.3,1.2){$x_1$};
\draw[ultra thick] (-0.5,1.5) -- (xyz cs:x=0.5,y=1.5);
\draw[ultra thick] (1,0) -- (xyz cs:x=1,y=1);
\node at (0.25,0.75){$M^b$};
   \end{tikzpicture} 
   \caption{Deformation in Example \ref{ex:ca} for $a< b <0$}
    \label{fig:1}
\end{figure}
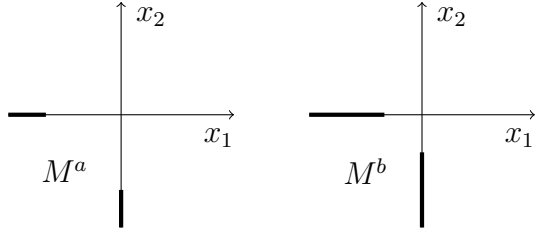

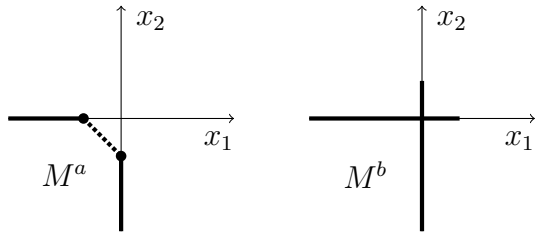
\begin{figure}
    \centering
\begin{tikzpicture}
\draw[->] (-3,0) -- (xyz cs:x=-3,y=3);
\node at (-2.6,2.8){$x_2$};
\draw[->] (-4.5,1.5) -- (xyz cs:x=-1.5,y=1.5);
\node at (-1.7,1.2){$x_1$};
\draw[ultra thick] (-4.5,1.5) -- (xyz cs:x=-3.5,y=1.5);
\draw[ultra thick] (-3,0) -- (xyz cs:x=-3,y=1);
\node at (-3.75,0.75){$M^a$};
\fill (-3.5,1.5) circle(2pt);
\fill (-3,1) circle(2pt);
\draw[ultra thick, densely dotted] (-3.5,1.5) -- (xyz cs:x=-3,y=1);

\draw[->] (1,0) -- (xyz cs:x=1,y=3);
\node at (1.4,2.8){$x_2$};
\draw[->] (-0.5,1.5) -- (xyz cs:x=2.5,y=1.5);
\node at (2.3,1.2){$x_1$};
\draw[ultra thick] (-0.5,1.5) -- (xyz cs:x=1.5,y=1.5);
\draw[ultra thick] (1,0) -- (xyz cs:x=1,y=2);
\node at (0.25,0.75){$M^b$};

   \end{tikzpicture} 
   \caption{Cell-attachment in Example \ref{ex:ca} for $a<0<b$}
    \label{fig:2}
\end{figure}

\begin{example}
  \label{ex:ca} 
  Consider the MPSC:
\begin{equation}
  \label{eq:ex1}
  \min \,\, f\left(x_1,x_2\right)=x_1+x_2 \quad \mbox{s.t.} \quad x_1 \cdot x_2 = 0. 
\end{equation}
We see that for all $a < b < 0$ the lower level sets $M^b$ and $M^a$ are homotopy-equivalent. More precisely, $M^b$ is a strong deformation retract of $M^a$, see Figure \ref{fig:1}. This is due to the fact that $M^b_a$ does not contain any $W$-stationary points. The situation changes dramatically if we pass the zero level corresponding to the W-stationary point $(0,0)$. For $a<0<b$ the lower level sets $M^b$ and $M^a$ are not homotopy-equivalent anymore. It is possible to describe the topological difference between $M^b$ and $M^a$ by means of the cell-attachment procedure. Namely, $M^b$ is homotopy-equivalent to $M^a$ with a one-dimensional cell attached along its boundary, see Figure \ref{fig:2}. \qed
\end{example}

In order to adequately describe the topological changes of MPSC lower level sets, we introduce nondegenerate W-stationary points along with their $W$-index. The latter will provide us with the dimension of the cell to be attached while passing the corresponding W-stationary level.

\begin{definition}[Nondegenerate W-stationarity]
   \label{def:nondeg_w_stationary}
   A W-stationary point $\bar x \in M$ with Lagrange multipliers $(\bar x, \bar \lambda,\bar \mu, \bar \sigma)$
   is called nondegenerate if the following conditions are satisfied:
   \begin{itemize}
        \item[] ND1: LICQ holds at $\bar x$,
        \item[] ND2: $\bar \mu_j > 0$ for all $j \in J_0(\bar x)$,
        \item[] ND3: $\bar \sigma_{1, m_\beta} \cdot \bar \sigma_{2, m_\beta} \ne 0$ for all $m_\beta \in \beta(\bar x)$,
        \item[] ND4: $D^2_{xx}L(\bar x, \bar \lambda,\bar \mu, \bar \sigma) \restriction_{T_{\bar x}M_0(\bar x)}$ is nonsingular.
   \end{itemize}
\end{definition}

Condition ND4 means that the matrix $V^T D^2L(\bar x) V$ is nonsingular, where $V$ is some matrix whose columns form a basis for the tangent space $T_{\bar x}M_0(\bar x)$ and $D^2_{xx}L(\bar x, \bar x, \bar \lambda,\bar \mu, \bar \sigma)$ is the Hessian of the Lagrange function with respect to $x$-variables.

\begin{definition}[W-index]
   \label{def:w_index}
   Let $\bar x \in M$ be a nondegenerate W-stationary point with unique Lagrange multipliers $(\bar x, \bar \lambda,\bar \mu, \bar \sigma)$.
   The number of negative eigenvalues of $D^2_{xx}L(\bar x, \bar \lambda,\bar \mu, \bar \sigma) \restriction_{T_{\bar x}M_0(\bar x)}$
   is called the quadratic index ($QI$) of $\bar x$. The number $\left|\beta(\bar x)\right|$ of non-zero bi-active pairs $(\bar \sigma_{1, m_\beta}, \bar \sigma_{2, m_\beta})$, $m_\beta \in \beta(\bar x)$ 
   is called the bi-active index ($BI$) of $\bar x$.
   The number $QI+BI$ is called the W-index of $\bar x$.
\end{definition}

Note that in absence of switching constraints, the W-index has only the $QI$-part and coincides
with the well-known quadratic index of a nondegenerate Karush-Kuhn-Tucker-point in nonlinear programming
or, equivalently, with the Morse index, cf. \cite{jongen:2000, kojima:1980, milnor:1963}.
Confer that for MPCC and MPVC the $BI$-part of the topologically relevant C-stationary and T-stationary points, respectively, counts the pairs of negative bi-active Lagrange multipliers, see \cite{jongen:2009} and \cite{dorsch:2012}. In MPSC, all pairs of bi-active Lagrange multipliers -- independently of their sign -- contribute to $BI$, which is an essentially new phenomenon here.    

The following Theorem \ref{theorem:Morse_lemma} describes the local structure of MPSC around a nondegenerate W-stationary point. 

\begin{theorem}[Morse Lemma for MPSC]
  \label{theorem:Morse_lemma}
  Suppose that $\bar x$ is a nondegenerate W-stationary point for MPSC with quadratic index $QI$, bi-active index $BI$ and
  W-index = $QI + BI$.
  Then, there exists a local $C^1$-coordinate system $\Psi: U \longrightarrow V$ of $\R^n$ around $\bar x$
   such that:
 \begin{equation}
    \label{eq:Morse_lemma}
    \begin{array}{rcl}
          f \circ \Psi^{-1} \left( 0_s, y_{s+1}, \ldots, y_n \right) &=& \displaystyle
    f(\bar x)+
    \Sum{j=1}{|J_0(\bar x)|} y_{j+s} \\
    &&\displaystyle +\Sum{m=1}{|\beta(\bar x)|} \left( y_{2m+q-1} + y_{2m+q} \right)
    +\Sum{k=1}{p} \pm y^2_{k+n-p},
    \end{array}
 \end{equation} 
where $y \in \{0_s\} \times \H^{|J_0(\bar x)|} \times \S^{|\beta(\bar x)|} \times \R^{p}$.
Moreover, in (\ref{eq:Morse_lemma}) there are exactly $QI$ negative squares. 
\end{theorem}
\proof
Without loss of generality, we may assume $f(\bar x) = 0$. Let $\Phi: U \longrightarrow V$ be a standard diffeomorphism
according to Definition (\ref{definition:standard_diffeo}). We put $\bar f=f \circ \Phi^{-1}$ on the set 
$\left( \{0_{s}\} \times \H^{|J_0(\bar x)|} \times \S^{|\beta(\bar x)|} \times \R^{p} \right) \cap V$. 
From now on we may assume $s=0$.
In view of Remark \ref{remark:deriv_in_new_coordinates} we have at the origin:
\begin{itemize}
\item[(i)]   $\displaystyle \frac{\partial \bar f}{\partial y_{j}} > 0$, $j \in J_0(\bar x)$,
\item[(ii)]  $\displaystyle \frac{\partial \bar f}{\partial y_{2m+q-1}} \cdot
             \frac{\partial \bar f}{\partial y_{2m+q}} \ne 0$, $m=1,\ldots |\beta(\bar x)|$,
\item[(iii)]  $\displaystyle \frac{\partial \bar f}{\partial y_{k+n-p}} = 0$, $k =1, \ldots, p$ and 
             $\left(
                   \displaystyle \frac{\partial^2 \bar f}{ \partial y_{k_1+n-p} \partial y_{k_2+n-p}}
             \right)_{1 \leq k_1, k_2 \leq p}$ is a nonsingular matrix with $QI$ negative eigenvalues.
\end{itemize}

From now on we denote $\bar f$ by $f$. Under the following coordinate transformations the set 
$\H^{|J_0(\bar x)|} \times \S^{|\beta(\bar x)|} \times \R^{p}$ will be 
transformed in itself (equivariance). As an abbreviation we put $y=\left( Y_{n-p}, Y^p \right)$,
where $Y_{n-p} = (y_1, \ldots, y_{n-p})$ and $Y^p = (y_{n-p+1}, \ldots, y_n)$. We write
\[
   f(Y_{n-p}, Y^p) = f(0, Y^p) + \int^1_0 \frac{d}{dt} f(t Y_{n-p}, Y^p) dt =
   f(0, Y^p) + \Sum{\ell=1}{n-p} y_\ell d_\ell(y), 
\]
where $d_\ell \in C^1$, $\ell=1,\ldots, n-p$.

In view of (iii) we may apply the Morse Lemma on the $C^2$-function $f(0, Y^p)$, cf. \cite[Theorem 2.8.2]{jongen:2000}, without affecting the coordinates $Y_{n-p}$. 
The corresponding coordinate transformation is of class $C^1$. Denoting the transformed functions again by $f$ and $d_\ell$, we obtain:
\[
   f(y) = \Sum{\ell=1}{n-p} y_\ell d_\ell(y) + \Sum{k=1}{p} \pm y^2_{k+n-p}.    
\]
Note that $d_\ell(0) = \displaystyle \frac{\partial f}{\partial y_\ell}(0)$, $\ell=1, \ldots, n-p$.
Recalling (i)-(ii), we put 
\begin{equation}
   \label{eq:last_trafo}
   \begin{array}{l}
         y_j |d_i(y)|, j=1,\ldots, |J_0(\bar x)|, \\
         y_m d_m(y), m=|J_0(\bar x)|+1,\ldots, n-p, \\  
         y_k, k = n-p+1, \ldots, n
   \end{array}
\end{equation}
as new local $C^1$-coordinates. Denoting the transformed function $f$ again by $f$ and, recalling the signs in (i)-(ii),
we obtain (\ref{eq:Morse_lemma}) by a straightforward application of the inverse function theorem. Here, the coordinate transformation $\Psi$ is understood as the composition of all previous ones. $\square$

As a by-product we elaborate on how the W-index can be used to characterize nondegenerate local minimizer of MPSC.

\begin{corollary}[Minimizers and W-index]
\label{lem:min-index}
Let $\bar x \in M$ be a nondegenerate W-stationary point. Then, $\bar x$ is a local
minimizer of MPSC if and only if its W-index vanishes.
\end{corollary}

\proof Let $\bar x$ be a nondegenerate W-stationary point for MPSC.
The application of Morse Lemma from Theorem \ref{theorem:Morse_lemma} says that there exist
neighborhoods $U$ and $V$ of $\bar x$ and $0$, respectively, and a local $C^1$-coordinate system $\Psi: U \rightarrow V$ of $\R^n$ around $\bar x$ such that (\ref{eq:Morse_lemma}) holds. Therefore, $\bar x$ is a local minimizer for MPSC if and only if $0$ is a local minimizer of $ f\circ \Psi^{-1}$ on the set $\{0_s\} \times \H^{|J_0(\bar x)|} \times \S^{|\beta(\bar x)|} \times \R^{p}$.
If the W-index vanishes, we have $\beta(\bar x)=\emptyset$ and $QI=0$, and (\ref{eq:Morse_lemma}) reads as
 \begin{equation}
    \label{eq:normalatMI0}
    f \circ \Psi^{-1} \left( 0_s, y_{s+1}, \ldots, y_n \right) =
    f(\bar x)+
    \Sum{j=1}{|J_0(\bar x)|} y_{j+s} 
    +\Sum{k=1}{p} y^2_{k+n-p},
 \end{equation}
where $y \in \{0_s\} \times \H^{|J_0(\bar x)|} \times \R^{p}$. Thus, $0$ is a local minimizer for (\ref{eq:normalatMI0}).
Vice versa, if $0$ is a local minimizer for (\ref{eq:Morse_lemma}), then obviously $\beta(\bar x)=\emptyset$ and $QI= 0$, hence, the W-index of $\bar x$ vanishes.
\qed

\begin{remark}[Nondegenerate minimizers] Corollary \ref{lem:min-index} says  that at a nondegenerate minimizer $\bar x \in M$ of MPSC with Lagrange multipliers $(\bar \lambda, \bar \mu, \bar \sigma)$ we have:
\begin{itemize}
    \item[1:]  LICQ holds at $\bar x$,
    \item[2:]  the Lagrange multipliers corresponding to active inequality constraints are positive, i.\,e. $\bar \mu_j >$ for all $j \in J_0(\bar x)$,
    \item[3:]  the bi-active switching constraints are absent, i.\,e. $\beta(\bar x)=\emptyset$,
    \item[4:]  the second-order sufficient condition (SOSC) is fulfilled, i.\,e. the matrix $D^2_{xx}L(\bar x, \bar \lambda,\bar \mu, \bar \sigma) \restriction_{T_{\bar x}M_0(\bar x)}$ is positive definite.
\end{itemize}
 Whereas conditions 1, 2, and 4 are standard as compared e.g. to nondegenerate minimizers of MPCC or MPVC, condition 3 is new and specific for MPSC. \qed
\end{remark}
 
The assumption of nondegeneracy is justified in the sense, that it generically holds at all W-stationary points of MPSC.

\begin{theorem}[Nondegeneracy is generic]
\label{thm:ND-gen}
Let $\mathcal{F}$ denote the subset of $C^2(\R^n,\R^{|I|}) \times C^2(\R^n,\R^{|J|}) \times C^2(\R^n,\R^k) \times C^2(\R^n,\R^k)$ consisting of those MPSC defining functions for which each 
W-stationary point is nondegenerate. Then, $\mathcal{F}$ is $C^2_s$-open and -dense.
\end{theorem}


\proof

Let us fix an index set $J_0 \subset \left\{1,\ldots,|J|\right\}$ of active inequality constraints, an index subset  $K_0 \subset J_0$ of these active inequality constraints, pairwise disjoint index sets $\alpha,\gamma, \beta \subset \{1,\ldots,k\}$ of switching constraints satisfying $\alpha \cup \gamma \cup \beta = \{1,\ldots,k\}$, an index subset  $\delta \subset \beta$ of be-active switching constraints,
and a number $r \in \mathbb{N}$ standing for the rank.
For this choice we consider the set  $M_{J_0,K_0,\alpha,\gamma,\beta, \delta,r}$ of $x \in \R^n$ such that the following conditions are satisfied:
\begin{itemize}
    \item[] (m1) $h_i(x)=0$ for all $i \in I$ and $g_{j}(x)=0$ for all  $j \in J$,
    \item[] (m2) $F_{1,m_\alpha}(x) =0$ for all $m_\alpha \in \alpha$, $F_{2,m_\gamma}(x) =0$ for all $m_\gamma \in \gamma$, $F_{1,m_\beta}(x) = F_{2,m_\beta}(x)=0$ for all $m_\beta \in \beta$,
    \item[] (m3) 
       \[\begin{array}{rcl}
              D f(x) &=&\displaystyle \Sum{i \in I}{} \lambda_i D h_i(x)+\Sum{j \in J_0\backslash K_0}{}  \mu_j D g_j( x) \\
        && \displaystyle +\Sum{m_\alpha \in \alpha}{} \sigma_{1,m_\alpha} D F_{1,m_\alpha}(x)+ \Sum{m_\gamma \in \gamma}{}
        \sigma_{2,m_\gamma} D F_{2,m_\gamma}(x)\\ 
        &&\displaystyle +\Sum{m_\beta \in \beta \backslash \delta}{} \left( \sigma_{1,m_\beta} D F_{1,m_\beta}(x)+
         \sigma_{2,m_\beta} D F_{2,m_\beta}(x) \right),
         \end{array} 
         \]
    \item[] (m4) the matrix $D^2_{xx} L( x,\lambda,\mu,\sigma)\restriction_{\mathcal{T}_{x} M(x)}$ has rank $r$.
\end{itemize}
Note that (m1) refers to equality and active inequality constraints, (m2) to switching constraints, while (m3) describes violation of ND2 and ND3. Furthermore, (m4) describes violation of ND4.
Now, it suffices to show that $M_{J_0,K_0,\alpha,\gamma,\beta, \delta,r}$ is generically empty whenever one of the sets $K_0$ or $\delta$ is nonempty or the rank $r$ in (m4) is not full, i.\,e. $r < \mbox{dim}\left(\mathcal{T}_{x} M(x)\right)$.
In fact, the available degrees of freedom
of the variables involved in each $M_{J_0,K_0,\alpha,\gamma,\beta, \delta,r}$ are equal to $n$. The loss of freedom caused by (m1) is $\left|I\right|+\left|J_0\right|$, and
the loss of freedom caused by (m2) is $\left|m_\alpha\right|+\left|m_\gamma\right|+2\left|m_\beta\right|$.
Due to Theorem \ref{thm:licq-gen}, LICQ holds generically at any feasible $x$, i.\,e. (ND1) is fulfilled. Suppose that the sets $K_0$ and $\delta$ are empty, then (m3) causes a loss of freedom of $n-\left|P\right|-\left|Q_0\right|-\left|m_\alpha\right|-\left|m_\gamma\right|-2\left|m_\beta\right|$. Hence, the total loss of freedom is $n$.
We conclude that a
further degeneracy, i.\,e. $K_0 \not = \emptyset$, $\delta \not = \emptyset$ or $r < \mbox{dim}\left(\mathcal{T}_{x} M(x)\right)$, would imply that the total available degrees of freedom $n$ are exceeded. By virtue of the jet
transversality theorem from \cite{jongen:2000}, generically the sets $M_{J_0,K_0,\alpha,\gamma,\beta, \delta,r}$ must be empty.
For the openness result, we argue in a standard way. Locally, W-stationarity can be written
via stable equations. Then, the implicit function theorem for Banach spaces can be applied to
follow W-stationary points with respect to (local) $C^2$-perturbations of defining functions. Finally,
a standard globalization procedure exploiting the specific properties of the strong $C^2_s$-topology can be used to construct a (global) $C^2_s$-neighborhood of problem data for which the nondegeneracy
property is stable, cf. \cite{jongen:2000}. \qed

Now, we are ready to state and prove deformation and cell-attachment theorems for MPSC, which constitute the core results of the paper. 

\begin{theorem}
 \label{theorem:main}
   Let $M^b_a$ be compact and suppose that LICQ is satisfied at all points $x\in M^b_a$.
\begin{itemize}
\item [(a)] {\bf (Deformation)} If $M^b_a$ does not
  contain any W-stationary point for MPSC, then $M^a$ is a
  strong deformation retract of $M^b$.

\item [(b)]  {\bf (Cell-attachment)} If $M^b_a$ contains exactly one W-stationary point for MPSC,
  say $\bar{x}$, and if $a<f(\bar{x})<b$ and the W-index of
  $\bar{x}$ is equal to $w$, then $M^b$ is homotopy-equivalent to
  $M^a$ with a $w$-cell attached.
  \end{itemize}
\end{theorem}

\proof (a) 
Due to LICQ at all $x\in M^b_a$ there exist real numbers
   $\lambda_i(x)$, $i \in I$,
   $\sigma_{1,m_\alpha}(x)$, $m_\alpha \in \alpha(x)$,
   $\sigma_{2,m_\gamma}(x)$, $m_\gamma \in \gamma(x)$,
   $\mu_j(x)$, $j \in J_0(x)$,
   $\sigma_{1, m_\beta}(x)$, $\sigma_{2, m_\beta}(x)$, $m_\beta \in \beta(x)$,
   $\nu_{l}(x)$, $l=1, \ldots , p$
     such that:
   \[
      \begin{array}{rcl}
            D f(x) &=&\displaystyle \Sum{i \in I}{} \lambda_i(x) D h_i(x)+
         \Sum{j \in J_0(x)}{} \mu_j(x) D g_j(x)
   \\
         &&\displaystyle+\Sum{m_\alpha \in \alpha(x)}{} \sigma_{1,m_\alpha}(x) D F_{1,m_\alpha}(x)+
         \Sum{m_\gamma \in \gamma(x)}{} \sigma_{2,m_\gamma}(x) D F_{2,m_\gamma}(x)
   \\
          &&\displaystyle +
         \Sum{m_\beta \in \beta(x)}{} \left( \sigma_{1,m_\beta}(x) D F_{1,m_\beta}(x)+
         \sigma_{2,m_\beta}(x) D F_{2,m_\beta}(x) \right)+
         \Sum{l=1}{p} \nu_{l}(x) \, \xi_{l},
      \end{array}
   \]
where vectors $\xi_{l}$, $l=1, \ldots,p$ are chosen as in Lemma \ref{lemma:feasible_set}.
We set:
\[
\begin{array}{rcl}
       A&=&\displaystyle \left\{\left. x \in M^b_a \,\right|\, \mbox{there exists } l \in \{1,\ldots, p\} \mbox{ with } \nu_{l}(x) \not = 0 \right\},
\\ \\
   B&=&\displaystyle \left\{\left. x \in M^b_a \,\right|\, \mbox{there exists } j \in J_0(x) \mbox{ with } \mu_j(x) <0 \right\}, 
\end{array}
\]
Since each $\bar x \in M^b_a$ is not W-stationary for MPSC, we get $\bar x \in A \cup B$.

The proof consists of a local argument and its globalization which are more or less standard here, cf. \cite[Theorem 3.2]{jongen:2009}. Let us briefly recall the details for the sake of completeness.

First, we show the {\bf local argument:} for each $\bar x \in M^b_a$ there exist an $(\R^n)$-neighborhood $U_{\bar x}$ of $\bar x$, ${t_{\bar x}} >0$ and a mapping
\[ 
\Psi^{\bar x}: \left\{ 
  \begin{array}{ccc}
   [0, t_{\bar x}) \times \left( M^b \cap U_{\bar x} \right) & \longrightarrow  & M  \\
   (t,x)                                        & \mapsto & \Psi^{\bar x}(t,x)
  \end{array}
         \right.  
\]
such that
\begin{itemize}
   \item[(i)]   $\Psi^{\bar x}(t, M^b \cap U_{\bar x}) \subset M^{b-t}$ for all $t \in [0, t_{\bar x})$, 
   \item[(ii)]  $\Psi^{\bar x}(t_1+t_2,\cdot)= \Psi^{\bar x}(t_1,\Psi^{\bar x}(t_2,\cdot))$ for all $t_1, t_2 \in [0, t_{\bar x})$ 
                with $t_1+t_2 \in [0, t_{\bar x})$,          
   \item[(iii)] $\Psi^{\bar x}(\cdot,\cdot)$ is 
                a $C^1$-flow corresponding to a $C^1$-vector field $F^{\bar x}$. 
\end{itemize}  
                 
Obviously, the level sets of $f$ are locally mapped onto the level sets of $f \circ \Phi^{-1}$, where $\Phi$ is a $C^1$-diffeomorphism according to Definition \ref{definition:local_coordinates}.
Applying the standard diffeomorphism $\Phi$ from Definition \ref{definition:standard_diffeo}, we consider
$f \circ \Phi^{-1}$ (denoted by $f$ again). Thus, we have $\bar x = 0$ and $f$ is given on the feasible set 
$\{0_s\} \times \H^{|J_0(\bar x)|} \times \S^{|\beta(\bar x)|} \times \R^{p}$.

{\bf Case a)} $\bar x \in A$\

Then, due to Remark \ref{remark:deriv_in_new_coordinates} there exists 
$l \in \{1,\ldots, p\}$ with $\displaystyle \frac{\partial f}{\partial x_l} (\bar x) \not = 0$. Define a local $C^1$-vector field $F^{\bar x}$ as follows:
\[
   F^{\bar x}(x_1, \ldots, x_l, \ldots, x_n)= 
   \left( 0, \ldots, -\frac{\partial f}{\partial x_l} (x) \cdot \left( \frac{\partial f}{\partial x_l} (x) \right)^{-2}, \ldots, 0 \right)^T.
\]
After respective inverse changes of local coordinates $F^{\bar x}$ induces the flow $\Psi^{\bar x}$, which
fits the local argument, see \cite[Theorem 2.7.6]{jongen:2000} for details.  

{\bf Case b)} $\bar x \in B$

Then, due to Remark \ref{remark:deriv_in_new_coordinates} there exists 
$j \in J_0(x)$ with $\displaystyle \frac{\partial f}{\partial x_j} (\bar x) < 0$. By means of a $C^1$-coordinate transformation in the $j$-th coordinate on $\H$, leaving the other coordinates unchanged, we obtain locally for $f$, see \cite[Theorem 3.2.26]{jongen:2000}: 
\[
  f(x_1, \ldots, x_j, \ldots, x_n) = - x_j + f(x_1, \ldots, \bar x_j, \ldots, x_n).
\]  
Define a local $C^1$-vector field $F^{\bar x}$ as follows:
\[
   F^{\bar x}(x_1, \ldots, x_j, \ldots, x_n)= \left( 0, \ldots, 1 , \ldots, 0 \right)^T.
\]
After respective inverse changes of local coordinates $F^{\bar x}$ induces the flow $\Psi^{\bar x}$, which
fits the local argument, see \cite[Theorem 3.3.25]{jongen:2000} for details.  
 
Second, the {\bf globalization} of the local argument will be performed. For that, consider the open covering $\{ U_{\bar x} \,|\, \bar x \in M^b_a \}$ of $M^b_a$. Since $M^b_a$ is compact, we get a finite 
open subcovering $\{ U_{\bar x_j} \,|\, \bar x_j \in M^b_a\}$ of $M^b_a$. Using a $C^\infty$-partition of unity $\{\phi_j\}$ subordinate to $\{ U_{\bar x_j} \,|\, \bar x_j \in M^b_a \}$
we define with $F^{\bar x_j}$, cf. Cases a),b), a $C^1$-vector field $ F= \Sum{j}{} \phi_j F^{\bar x_j}$. 
The last induces a flow $\Psi$ on $\{ U_{\bar x_j} \,|\, \bar x_j \in M^b_a \}$, see \cite[Theorem 3.3.14]{jongen:2000} for details.
Moving along the local pieces of the trajectories $\Psi(\cdot,x)$, $x \in M^b_a$ reduces the level of $f$ at least by a positive real number
\[
 \frac{\min\{t_{x_j}\,|\, x_j \in M^b_a \}}{2}.
\]
Thus, we obtain for $x \in M^b_a$ a unique $t_a(x) > 0$ with $\Psi(t_a(x),x) \in M^a$. 
It is not hard (but technical) to realize that $t_a: x \longrightarrow t_a(x) $ is Lipschitz continuous.
Finally, we define $r:[0,1] \times M^b \longrightarrow M^b$ as follows:
\[ r(\tau, x)= \left\{ 
  \begin{array}{lll}
   x                     & \mbox{for } x \in M^a, & \tau \in [0,1]\\
   \Psi(\tau t_a(x),x)   & \mbox{for } x \in M^b_a, & \tau \in [0,1].
  \end{array}
         \right.  
\]
The mapping $r$ provides that $M^a$ is a strong deformation retract of $M^b$.

(b) Part (a) allows deformations up to an arbitrarily small neighborhood of the W-stationary point $\bar x$. In such a neighborhood, we may assume without loss of generality that $\bar x=0$ and $f$ has the following form as from Theorem \ref{theorem:Morse_lemma} (for simplicity we omit the trivial part by assuming $s=0$):
\begin{equation}
    \label{eq:n1}
 f (x) =
    f(\bar x)+
    \Sum{j=1}{|J_0(\bar x)|} x_{j} 
    +\Sum{m=1}{|\beta(\bar x)|} \left( x_{2m+q-1} + x_{2m+q} \right)
    +\Sum{k=1}{p} \pm x^2_{k+n-p},
\end{equation}
where $x \in \H^{|J_0(\bar x)|} \times \S^{|\beta(\bar x)|} \times \R^{p}$, and the number of negative squares in (\ref{eq:n1}) equals $QI$.
In terms of \cite{goresky:1988} the set 
$\H^{|J_0(\bar x)|} \times \S^{|\beta(\bar x)|} \times \R^{p}$ can be interpreted as the product of the tangential part $\R^p$ and the normal part $\H^{|J_0(\bar x)|} \times \S^{|\beta(\bar x)|}$. The cell-attachment along the tangential part is standard.  Analogously to the unconstrained case, a $QI$-dimensional cell has to be attached on $\R^p$. The cell-attachment along the normal part is more involved. First, we emphasize that the linear terms $x_j$, $j \in J_0(\bar x)$, in (\ref{eq:n1}) do not contribute to the dimension of the cell to be attached. In fact, with respect to lower level sets, the $1$-dimensional  constrained singularity $x$ on $\H$ plays the same role as the unconstrained singularity $x^2$. In this sense the linear terms on $\H^{|J_0(\bar x)|}$ can be neglected for cell-attachment. Second, the linear terms $x_{2m+q-1} + x_{2m+q}$, $m \in \beta(\bar x)$, in (\ref{eq:n1}) contribute each with the one-dimensional cell due to Example \ref{ex:ca}. Hence, the linear terms on $\S^{|\beta(\bar x)|}$ are responsible for the attachment of a $BI$-dimensional cell -- recall that $BI=|\beta(\bar x)|$.  Finally, we apply Theorem 3.7 from Part I in \cite{goresky:1988}, which says that the local Morse data is the product of tangential and normal Morse data. Hence, the dimensions of the attached cells add together. Here, a $(QI+BI)$-dimensional cell has to be attached. The latter corresponds to the W-index of the W-stationary point $\bar x$. \qed

Let us present a global interpretation of Theorem \ref{theorem:main} by showing a mountain pass theorem for MPSC. 

\begin{theorem}[Mountain pass]
\label{thm:mrel}
Let the MPSC feasible set $M$ be compact and connected, and all W-stationary points of MPSC be nondegenerate. 
Then, it holds:
\begin{equation}
    \label{eq:mr}
    r_s \geq r-1,
\end{equation}
where $r$ is the number of local minimizers of MPSC and $r_s$ is the number of W-stationary saddle points with W-index equal to one.
\end{theorem}

\proof 
We assume without loss of generality that the objective function $f$ has pairwise different values at all W-stationarity points of MPSC. If it is not the case, we may enforce this property by sufficiently small perturbations of the objective function.  
Due to the openness part in Theorem \ref{thm:ND-gen}, all W-stationarity points of such a perturbed MPSC remain nondegenerate. Moreover, the formula (\ref{eq:mr}) is still valid since it does not depend on the functional values of $f$.

Further, let $q_a$ denote the number of connected components of the lower level set $M^a$.
We focus on how $q_a$ changes as $a \in \R$ increases. Due to Theorem \ref{theorem:main}(a), $q_a$ can change only if passing through a value corresponding to a W-stationary point $\bar x$, i.\,e. $a=f\left(\bar x\right)$. 
In fact, Theorem \ref{theorem:main}(a) allows homotopic deformations of lower level sets up to an arbitrarily small neighborhood of the W-stationary point $\bar x$.
Then, we have to estimate the difference between $q_a$ and $q_{a-\varepsilon}$, where $\varepsilon > 0$ is arbitrarily, but sufficiently small, and $a=f\left(\bar x\right)$. This is done by a {\bf local argument}. For that, let the W-index of $\bar x$ be $QI+BI$. We use Theorem \ref{theorem:main}(b) which says that $M^{a}$ is homotopy-equivalent to $M^{a-\varepsilon}$ with a $(QI+BI)$-dimensional cell attached.
Let us distinguish the following cases:
\begin{itemize}
    \item[1)] $\bar x$ is a local minimizer with vanishing W-index, cf. Corollary \ref{lem:min-index}, i.\,e. $QI=BI=0$. Then, we attach to $M^{a-\varepsilon}$ the cell of dimension zero. Consequently, a new connected component is created, and it holds: 
\[
   q_a = q_{a-\varepsilon} + 1.
\]
    \item[2)] $\bar x$ is a saddle point with W-index equal to one, i.\,e. either $QI=1$, $BI=0$ or $QI=0$, $BI=1$ holds.
    Then, we attach to $M^{a-\varepsilon}$ the cell of dimension one. Consequently, at most one connected component disappears, and it holds:
\[
   q_{a-\varepsilon} -1 \leq q_a \leq q_{a-\varepsilon}.
\]
 The case $QI=1$, $BI=0$ is well known from nonlinear programming, see e.\,g. \cite{jongen:2000}.
 The case $QI=0$, $BI=1$ is new and characteristic for MPSC, cf. Example \ref{ex:srel}.
    \item[4)] $\bar x$ is W-stationary with W-index greater than one, i.\,e. $QI+BI > 1$.
    The boundary of the to be attached $(QI+BI)$-dimensional cell is thus connected. Consequently, the number of connected components of $M^{a}$ remains unchanged, and it holds:
\[
   q_a = q_{a-\varepsilon}.
\]
\end{itemize}

Now, we proceed with the {\bf global argument}. Compactness of $M$ implies that there exists $c \in \R$ such that $M^c$ is empty, thus, $q_c =0$. Additionally, there exists $d \in \R$ such that $M^d$ is connected and contains all W-stationary points, thus, $q_d=1$. Since nondegenerate W-stationary points are in particular isolated, we conclude that there must be finitely many of them again due to the compactness of $M$. 
Let us now increase the level $a$ from $c$ to $d$ and describe how the number $q_a$ of connected components of the lower level sets $M^a$ changes.
It follows from the local argument that $r$ new connected components are created, where $r$ is the number of local minimizers for MPSC. Let $q_s$ denote the actual number of disappearing connected components if passing the levels corresponding to W-stationary saddle points.
The local argument provides that at most $r_s$ connected components might disappear while doing so, i.\,e.
\[
  q_s \leq r.
\]
Altogether, we have:
\[
  r - r_s \leq r - q_s = q_d-q_c.
\]
By recalling that $q_d=1$ and $q_c=0$, we get (\ref{eq:mr}). \qed

Finally, we link our results to the numerical investigations of MPSC undertaken in \cite{kanzow:2019}. There, the important role of W-stationary points is highlighted if applying the Scholtes' relaxation to switching constraints. In \cite[Theorem 3.2]{kanzow:2019} it is proven that KKT-points of such relaxation converge to a W-stationary point of MPSC under LICQ. Moreover, as shown in \cite[Example 3.3]{kanzow:2019} stronger results
cannot be expect in general. Let us illustrate how the global structure of MPSC is inherited by the Scholtes' relaxation.   

\begin{example}[Scholtes' relaxation, cf. \cite{kanzow:2019}]
  \label{ex:srel} 
  Consider the MPSC:
\begin{equation}
  \label{eq:ex3}
  \min \,\, f\left(x_1,x_2\right)=\left(x_1-1\right)^2+\left(x_2-1\right)^2 \quad \mbox{s.t.} \quad x_1 \cdot x_2 = 0. 
\end{equation}
It is easy to see that the feasible point $x^1=(0,0)$ is W-stationary. Moreover, it is nondegenerate with quadratic index $QI=0$. For its bi-active index we have $BI=1$.
This means that $x^1$ is a saddle point with W-index equal to one which connects two minimizers $x^2=(1,0)$ and $x^3=(0,1)$. Let us consider the Scholtes' relaxation for (\ref{ex:srel}) with sufficiently small $t >0$:
\begin{equation}
  \label{eq:ex3-rel}
  \min \,\, f\left(x_1,x_2\right)=\left(x_1-1\right)^2+\left(x_2-1\right)^2 \quad \mbox{s.t.} \quad -t \leq x_1 \cdot x_2 \leq t. 
\end{equation}
It is easily seen that the KKT-points of (\ref{eq:ex3-rel}) are 
  \[
      x^1_t = \left( \sqrt{t}, \sqrt{t} \right),
  \]
  \[
      x^2_t = 
      \left( \displaystyle \frac{1+\sqrt{1 - 4 t}}{2}, \displaystyle \frac{1-\sqrt{1 - 4 t}}{2} \right),
  \]
  \[
      x^3_t = 
      \left( \displaystyle \frac{1-\sqrt{1 - 4 t}}{2}, \displaystyle \frac{1+\sqrt{1 - 4 t}}{2} \right).
  \]
  Obviously, $x^1_t \longrightarrow x_1$, $x^2_t \longrightarrow x_2$, and $x^3_t  \longrightarrow x_3$ as $t \longrightarrow 0$.
  Moreover, $x^2_t$ and $x^3_t$ are minimizers for 
  (\ref{eq:ex3-rel}) with zero quadratic index, 
  the quadratic index of $x^1_t$ is equal to one.
  By the smoothing procedure, the W-stationary point $x^1$ with W-index equal to one
  corresponds to the KKT-point $x^1_t$ with quadratic index equal to one likewise. 
  We conclude that the smoothed version fully preserves the global topological structure of MPSC.
\end{example}

\section{Characterization of strong stability}

We would like to address the question whether the requirements ND1-ND4 for W-stationary points are too restrictive. For that, let us examine the gap between nondegeneracy and strong stability of W-stationary points. The concept of strong stability is defined by means of an appropriate semi-norm. To this aim let be $\bar x \in \R^n$, $r>0$. For defining functions $\left(f, h, g, F_1, F_2\right)$ from (\ref{eq:math-program}) the seminorm 
$\left\|\left(f, h, g, F_1, F_2 \right)\right\|^{C^2}_{B(\bar x, r)}$
is defined to be the maximum modulus of the function values and partial derivatives up to order two
of $f, h, g, F_1, F_2$.

\begin{definition}[Strong Stability, cf. \cite{kojima:1980}]
  \label{def:strongly_stable_point}
  A W-stationary point $\bar x \in M$, resp. a W-stationary pair 
$(\bar x, \bar \lambda, \bar \mu, \bar \sigma)$, 
 for $MPSC[f,g,h,F_1,F_2]$ is called ($C^2$)-strongly stable 
if for some $r > 0$ and each $\varepsilon \in (0,r]$ there exists $\delta = \delta(\varepsilon) > 0$
such that whenever $\left( \widetilde f, \widetilde h, \widetilde g, \widetilde F_1, \widetilde F_2 \right) \in C^2$ and 
\[
\left\| \left( f- \widetilde f, h - \widetilde h, g - \widetilde g, F_1 - \widetilde F_1, F_2 - \widetilde F_2 \right) \right\|^{C^2}_{B(\bar x, r)} \leq \delta, 
\]
 $B(\bar x, \varepsilon)$, resp.
$B\left(\left(\bar x, \bar \lambda, \bar \mu, \bar \sigma \right), \varepsilon \right)$,
 contains a W-stationary point 
$\widetilde x$, resp. a W-stationary pair 
$\left(\widetilde x, \widetilde \lambda, \widetilde \mu, \widetilde \sigma \right)$,
for MPSC$\left[ \widetilde f, \widetilde h, \widetilde g, \widetilde F_1, \widetilde F_2  \right]$
which is unique in $B(\bar x, r)$, resp. 
unique in $B \left(\left(\widetilde x, \widetilde \lambda, \widetilde \mu, \widetilde \sigma \right), r \right)$.
\end{definition}

The following Lemma \ref{lemma:point_pair} establishes the connection between both definitions just introduced, cf. \cite{klatte:1990} for the case of nonlinear programming.

\begin{lemma}[W-stationary points and pairs]
  \label{lemma:point_pair}
The following assertions are equivalent:
\begin{itemize}
 \item[(a)] $\bar x$ is a strongly stable W-stationary point for MPSC which satisfies LICQ, and 
            $(\bar \lambda, \bar \mu, \bar \sigma)$ is the associated Lagrange multiplier vector.
 \item[(b)] $(\bar x, \bar \lambda, \bar \mu, \bar \sigma)$ is a strongly stable W-stationary pair for MPSC.
\end{itemize}
\end{lemma}

\proof 
(a) $\PPfeil$ (b) LICQ remains valid under small perturbations of the defining functions.
Hence, the corresponding Lagrange multipliers are unique. 
Remark \ref{remark:deriv_in_new_coordinates} 
provides the continuity of Lagrange multipliers with respect to perturbations under consideration.

(b) $\PPfeil$ (a) The nontrivial part is to prove that LICQ holds at $\bar x$. 
The proof goes along the lines of \cite[Theorem 2.3]{klatte:1990}. 
To stress the new aspects here we assume that there are only bi-active switching constraints, i.\,e. $I = \emptyset$, $J=\emptyset$, $\alpha(\bar x) = \emptyset$ and $\gamma(\bar x) = \emptyset$.
Let $(\bar x, \bar \sigma)$ be a strongly stable W-stationary pair for MPSC and 
let LICQ be not fulfilled at $\bar x$. Then, there exist
real numbers
   $\delta_{1,m}, \delta_{1,m}$, $m \in \beta (\bar x)$ (not all vanishing)
such that:
      \begin{equation}
          \label{eq:beweis_lemma}
              \Sum{m \in \beta(\bar x)}{} \left( \delta_{1,m} D F_{1,m}(\bar x)+
          \delta_{2,m} D F_{2,m}(\bar x) \right) = 0.
      \end{equation} 
We define
\[
  c = \Sum{m \in \beta(\bar x)}{} \left( D F_{1,m}(\bar x)+
          D F_{2,m}(\bar x) \right).
\]
For $\varepsilon > 0$ let
\[
   \sigma_{1,m}(\varepsilon) = \bar \sigma_{1,m} + \varepsilon,  
   \sigma_{2,m}(\varepsilon) = \bar \sigma_{2,m} + \varepsilon \mbox{ for all }
   m \in m_\beta(\bar x).
\]
Putting $\varphi(x) = c \cdot x$ we obtain that $(\bar x, \sigma(\varepsilon))$ is a W-stationary pair for 
MPSC$\left[ f+\varepsilon \cdot \varphi, F_1, F_2  \right]$.
Moreover, due to the strong stability of $(\bar x, \bar \sigma)$ for MPSC$\left[ f, F_1, F_2  \right]$
we claim that for each sufficiently small $\varepsilon > 0$ the W-stationary pair $(\bar x, \sigma(\varepsilon))$
is unique for MPSC$\left[ f+\varepsilon \cdot \varphi, F_1, F_2  \right]$ in some neighborhood $U$ of $(\bar x, \bar \sigma)$. 
However, (\ref{eq:beweis_lemma}) and $\sigma_{i,m}(\varepsilon) \not = 0 $ for $m \in m_\beta(\bar x)$, $i=1,2$,
ensure that for any sufficiently small real number $t$, the pair $(\bar x, v(\varepsilon, \delta, t))$ with
\[
   v_{1,m}(\varepsilon, \delta, t) =  \sigma_{1,m}(\varepsilon) + \delta_{1,m} t,
   v_{2,m}(\varepsilon, \delta, t)  = \sigma_{2,m}(\varepsilon) + \delta_{2,m} t \mbox{ for all }
   m \in m_\beta(\bar x)
\]
belongs to $U$ and is a W-stationary pair for MPSC$\left[ f+\varepsilon \cdot \varphi, F_1, F_2  \right]$.
Hence, necessarily $\delta = 0$, and LICQ is shown.
\qed

Now we give two guiding examples for instability which may occur at W-stationary points.

\begin{example}[Instability I]
\label{ex:inst1}
Consider the MPSC:
\begin{equation}
  \label{eq:exinst1}
  \min \,\, \pm x_1^2 \pm x_2^2 \quad \mbox{s.t.} \quad x_1 \cdot x_2 = 0. 
\end{equation}
Obviously, $(0,0)$ is the unique W-stationary point for (\ref{eq:exinst1}) with both vanishing bi-active Lagrange multipliers.
Consider the following perturbation of (\ref{eq:exinst1}) with respect to parameter $t>0$:
\begin{equation}
  \label{eq:ex1_pert}
  \min \,\, \pm \left(x_1-t\right)^2 \pm \left(x_2-t\right)^2 \quad \mbox{s.t.} \quad x_1 \cdot x_2 = 0.
\end{equation}
It is easy to see that $(0,0)$, $(0,t)$ and $(t,0)$ are W-stationary points for (\ref{eq:ex1_pert}). 
It means that $(0,0)$ is not a strongly stable W-stationary point for (\ref{eq:exinst1}). \qed
\end{example}

\begin{example}[Instability II]
\label{ex:inst2}
Consider the MPSC:
\begin{equation}
  \label{eq:exinst2}
  \min \,\, x_1 \pm x_2^2 \quad \mbox{s.t.} \quad x_1 \cdot x_2 = 0. 
\end{equation}
Obviously, $(0,0)$ is the unique W-stationary point for (\ref{eq:exinst1}) with exactly one vanishing bi-active Lagrange multiplier.
Consider the following perturbation of (\ref{eq:exinst2}) with respect to parameter $t>0$:
\begin{equation}
  \label{eq:ex2_pert}
  \min \,\, x_1 \pm \left(x_2-t\right)^2 \quad \mbox{s.t.} \quad x_1 \cdot x_2 = 0.
\end{equation}
It is easy to see that $(0,0)$ and $(0,t)$ are W-stationary points for (\ref{eq:ex2_pert}). 
It means that $(0,0)$ is not a strongly stable W-stationary point for (\ref{eq:exinst2}). \qed
\end{example}

For characterizing strong stability, we shall use some auxiliary objects associated with a W-stationary point $\bar x \in M$ and its multipliers $(\bar \lambda,\bar \mu, \bar \gamma)$. Let the index set of positive multipliers corresponding to the inequality constraints is given by
\[
J_+(\bar x) = \left\{ q \in J_0(\bar x) \,\left|\, \bar \mu_j > 0\right.\right\}.
\]
For $J_+(\bar x) \subset J_*\subset J_0(\bar x)$ we set 
\[
  \begin{array}{ll}
      M_*(\bar x)= \{ x \in \R^n \,\,| & h_i(x) =0, i \in I,
    g_j(x) =0, j \in J_*,  \\ &
       F_{1,m_\alpha}(x) = 0, m_\alpha \in \alpha(\bar x), F_{2,m_\gamma}(x) = 0, m_\gamma \in \gamma(\bar x), \\ &
    F_{1,m_\beta}(x) = 0, F_{2,m_\beta}(x) = 0, m_\beta \in \beta(\bar x)
       \}.
     \end{array}
\]
Obviously, $M_*(\bar x) \subset M$ and, in the case where LICQ holds at $\bar x$, $M_*(\bar x)$ is locally a
$C^2$-manifold of dimension $p + \left|J_0(\bar x)\right| - \left|J_*(\bar x)\right|$. The tangent space of $M_*(\bar x)$ at $\bar x$ is
\[
    \begin{array}{ll}
      T_{\bar x}M_*(\bar x) = \{ \xi \in \R^n \,\,| &
        D h_i(\bar x) \, \xi=0, i \in I, D g_j(\bar x) \, \xi =0, j \in J_* \\ &
        D F_{1,m_\alpha}(\bar x) \, \xi = 0, m_\alpha \in \alpha(\bar x), \\ &
        D F_{2,m_\gamma}(\bar x) \, \xi= 0, m_\gamma \in \gamma(\bar x), 
     \\ &
    D F_{1,m_\beta}(\bar x) \, \xi = 0, D F_{2,m_\beta}(\bar x) \, \xi = 0, m_\beta \in \beta(\bar x)
         \}.
     \end{array}
\]

\begin{theorem}[Characterization of strong stability]
\label{thm:char-ss}
Let $\bar x \in M$ be a W-stationary point of MPSC satisfying LICQ with the unique Lagrange multipliers $(\bar \lambda, \bar \mu, \bar \sigma)$. Then, $\bar x$ is strongly stable if and only if 
\begin{itemize}
    \item[(i)] it fulfils ND3, i.\,e. $\bar \sigma_{1, m_\beta} \cdot \bar \sigma_{2, m_\beta} \ne 0$ for all $m_\beta \in \beta(\bar x)$, and
    \item[(ii)]  the matrices $D^2_{xx} L(\bar x, \bar \lambda, \bar \mu, \bar \sigma)\restriction_{\mathcal{T}_{\bar x} M_*(\bar x)}$ are nonsingular with the same determinant sign for all index subsets $J_+(\bar x) \subset J_* \subset J_0(\bar x)$. 
\end{itemize}
\end{theorem}

\proof In virtue of LICQ at $\bar x$, Lemma \ref{lemma:point_pair} allows us to deal equivalently with the strong stability of
the W-stationary pair $(\bar x, \bar \lambda, \bar \mu, \bar \sigma)$.

{\bf Case 1: ND3 holds}

We consider the following mapping $\TT:\R^{n+|I|+|J|+2k} \longrightarrow \R^{n+|I|+|J|+2k}$
locally at its zero $(\bar x, \bar \lambda, \bar \mu, \bar \sigma)$:
\[
   \TT(x, \lambda, \mu, \sigma)=
      \left(
          \begin{array}{c}
               D_x L(x, \lambda, \mu, \sigma) \\
               h(x) \\
               \min \left\{ \mu, g(x)\right\} \\
               F_{1, \alpha(\bar x)} (x)\\
               F_{2, \gamma(\bar x)} (x)\\
               F_{1, \beta(\bar x)} (x)\\
               F_{2, \beta(\bar x)} (x)\\
          \end{array}
      \right).
\]
Note that W-stationary pairs for MPSC - in a sufficiently small neighborhood of $(\bar x, \bar \lambda, \bar \mu, \bar \sigma)$ - are precisely the zeros of $\TT$. Moreover, characterization of strong stability for KKT-points in the case of nonlinear programming can be applied here. Indeed, as in \cite[Theorem 4.3]{klatte:1990}, the strong stability for $(\bar x, \bar \lambda, \bar \mu, \bar \sigma)$
can be characterized by the fact that all matrices in the Clarke's subdifferential 
$\partial \TT(\bar x, \bar \lambda, \bar \mu, \bar \sigma)$ are nonsingular.
The latter can be equivalently rewritten as condition (ii), cf. also \cite{kummer:1999} for the case of nonlinear programming.

{\bf Case 2: ND3 fails}

Let $\Phi: U \longrightarrow V$ be the standard diffeomorphism from Defnition \ref{definition:standard_diffeo}. 
We put $\bar f=f \circ \Phi^{-1}$ on the set 
$\left( \{0_{s}\} \times \H^{|J_0(\bar x)|} \times \S^{|\beta(\bar x)|} \times \R^{p} \right) \cap V$.
From now on we may assume $s=0$. In view of Remark \ref{remark:deriv_in_new_coordinates} and failure of ND3 we have at the origin:
\begin{itemize}
\item[(i)]   $\displaystyle \frac{\partial \bar f}{\partial y_{j}} \geq 0, j \in J_0(\bar x),$
\item[(ii)]  $\displaystyle \frac{\partial \bar f}{\partial y_{2m+q-1}} \cdot
             \frac{\partial \bar f}{\partial y_{2m+q}} = 0$ for at least one $m=1,\ldots |\beta(\bar x)|,$
\item[(iii)]  $\displaystyle \frac{\partial \bar f}{\partial y_{k+n-p}} = 0$, $k =1, \ldots, p$.
\end{itemize}
Without loss of generality, we assume $m=1$ in (ii).
From now on we denote $\bar f$ again by $f$.
In what follows, we successively perform arbitrarily small perturbations of $f$ such that the origin remains a W-stationary point on $\H^{|J_0(\bar x)|} \times \S^{|\beta(\bar x)|} \times \R^{p}$. 

1) As a stabilization step we add to $f$ an arbitrarily small linear-quadratic term 
\[
   \Sum{j=1}{|J_0(\bar x)|} c_{j} \cdot y_{j} 
    +\Sum{m=2}{|\beta(\bar x)|} \left( c_{2m+q-1} \cdot y_{2m+q-1} + c_{2m+q} \cdot y_{2m+q} \right)
    +\Sum{k=1}{p} c_{k+n-p} \cdot y^2_{k+n-p},
\]
such that it holds for the perturbed function (denoted again by $f$):
\begin{itemize}
\item[(i)]   $\displaystyle \frac{\partial f}{\partial y_{j}} > 0$, $j \in J_0(\bar x)$,
\item[(ii)]  $\displaystyle \frac{\partial f}{\partial y_{2m+q-1}} \cdot
             \frac{\partial f}{\partial y_{2m+q}} \ne 0$, $m=2,\ldots, |\beta(\bar x)|$,
\item[(iii)] $\displaystyle \frac{\partial f}{\partial y_{k+n-p}} = 0$, $k =1, \ldots, p$ and
         $\left(
                   \displaystyle \frac{\partial^2 f}{ \partial y_{k_1+n-p} \partial y_{k_2+n-p}}
             \right)_{1 \leq k_1, k_2 \leq p}$ is a nonsingular matrix.             
\item[(iv)]  $\displaystyle \frac{\partial f}{\partial y_{q+1}} \cdot \frac{\partial f}{\partial y_{q+2}} = 0$.
\end{itemize}

2) We approximate $f$ by means of a $C^\infty$-function in a small $C^2$-neighborhood of $f$ 
such leaving its value, first and second order derivatives at the origin invariant.
This can be done since $C^\infty$-functions lie $C^2$-dense in $C^2$-functions.
We denote the latter $C^\infty$-approximation again by $f$.

Due to the stabilization step 1) and step 2), we can restrict our considerations to the following situation: 
\[
  f \in C^\infty \left(\R^2,\R \right), 0 \mbox{ is a W-stationary point for } f{\restriction_{\S}} \mbox{ and }
    \displaystyle \frac{\partial f}{\partial x_1}(0) \cdot \frac{\partial f}{\partial x_2}(0) = 0.
\]

{\bf Case I:  $\displaystyle \frac{\partial f}{\partial x_1}(0) = \frac{\partial f}{\partial x_2}(0) = 0$.} 

We can write
\[
  f\left(x_1,x_2\right) = g_{1, 1}\left(x_1,x_2\right) x^2_1 + 2 g_{1, 2}\left(x_1,x_2\right) x_1 x_2 + g_{2, 2}\left(x_1,x_2\right) x_2^2
\]
with $g_{1, 1}, g_{1, 2}, g_{2, 2} \in C^\infty \left( \R^2, \R \right)$.
Adding to $f$ an arbitrarily small quadratic term $ax_1^2+bx_2^2$, $a, b \in \R$ we can ensure that
$g_{1, 1}(0,0) \not = 0$ and $g_{2, 2}(0,0) \not = 0$.
Hence, 
\[
  \Psi\left(x_1,x_2\right)=\left(
               \begin{array}{c}
                   \displaystyle x_1 \cdot \sqrt{|{g_{1,1}\left(x_1,x_2\right)}|} \\
                   \displaystyle x_2 \cdot \sqrt{|{g_{2,2}\left(x_1,x_2\right)}|} \\
               \end{array}
             \right)
\]
is a local $C^\infty$-diffeomorphism leaving $\S$ invariant.
In new local coordinates induced by $\Psi$ we obtain:
\[
  f\left(x_1,x_2\right) = \varepsilon_1 x_1^2 + G\left(x_1,x_2\right) x_1 x_2 + \varepsilon_2 x_2^2,
\]
where $\varepsilon_1 = \mbox{sign}\left( g_{1, 1}(0,0) \right)$, $\varepsilon_2 = \mbox{sign}\left( g_{2, 2}(0,0) \right)$.
Since $G\left(x_1,x_2\right) x_1 x_2 = 0$ on $\S$, we can perturb $f$ by means of a real parameter 
as in Example \ref{ex:inst1} to get a bifurcation of $0$ as a W-stationary point.

{\bf Case II:  $\displaystyle \frac{\partial f}{\partial x_1}(0) \ne 0$, $\displaystyle \frac{\partial f}{\partial x_2}(0) = 0$.}

We can write
\[
  f\left(x_1,x_2\right) = g_{1, 1}\left(x_1,x_2\right) x_1 + g_{2, 2}\left(x_1,x_2\right) x_2^2
\]
with $g_{1, 1}, g_{2, 2} \in C^\infty \left( \R^2, \R \right)$ and $g_1(0,0) \ne 0$.
Adding to $f$ an arbitrarily small quadratic term $bx_2^2$, $b \in \R$ we can also ensure that $g_{2, 2}(0,0) \not = 0$. Hence, 
\[
  \Psi\left(x_1,x_2\right)=\left(
               \begin{array}{l}
                   \displaystyle x_1 \cdot g_{1,1}\left(x_1,x_2\right) \\
                   \displaystyle x_2 \cdot \sqrt{|{g_{2,2}\left(x_1,x_2\right)}|} \\
               \end{array}
             \right)
\]
is a local $C^\infty$-diffeomorphism leaving $\S$ invariant.
In new local coordinates induced by $\Psi$ we obtain:
\[
  f\left(x_1,x_2\right) = x_1 + \varepsilon x_2^2,
\]
where $\varepsilon = \mbox{sign}\left( g_{2, 2}(0,0) \right)$.
We can perturb $f$ by means of a real parameter 
as in Example \ref{ex:inst2} to get a bifurcation of $0$ as a W-stationary point. 

The case $\displaystyle \frac{\partial f}{\partial x_1}(0) = 0$, $\displaystyle \frac{\partial f}{\partial x_2}(0) \ne 0$ can be treated analogously.

Finally, performing all perturbations described above we ensure that $0$ is not a strongly stable W-stationary point. \qed

Next result is a direct consequences of Theorem \ref{thm:char-ss}. 

\begin{corollary}
 \label{corollary:ND-W-stationary}
Any nondegenerate W-stationary point of MPSC is strongly stable.
\end{corollary} 

Additionally, it is important to note that ND3 is necessary for strong stability of W-stationary points due to Theorem \ref{thm:char-ss}, i.\,e. both bi-active Lagrange multipliers cannot vanish. This new issue is in strong contrast e.\,g. with the characterization of strong stability of C-stationary points in MPCC, see \cite{jongen:2012}. In MPCC, strong stability includes cases where one of the bi-active Lagrange multipliers may vanish. In other words, there is an essential discrepancy between nondegeneracy and strong stability of C-stationary points in MPCC. In MPSC, the situation is rather different: in absence of active inequality constraints, nondegeneracy of W-stationary points is equivalent to their strong stability.

\begin{corollary}
 \label{corollary:ND-W-stationary1}
Let inequality constraints be inactive at a W-stationary point $\bar x \in M$, i.\,e. $J_0(\bar x)=\emptyset$. 
Then, $\bar x$ is nondegenerate if and only if it is strongly stable and satisfies LICQ.
\end{corollary}

\end{document}